\documentclass[a4paper,12pt]{amsart}

\usepackage[latin1]{inputenc}
\usepackage{amsxtra}
\usepackage{amssymb}
\usepackage{bbm}
\usepackage{program}
\sfvariables
\usepackage[all]{xypic}
\xyoption{dvips}
\usepackage{longtable}
\usepackage{array}
\usepackage{epsfig}

\theoremstyle{definition}
\newtheorem{ntn}{Notation}[section]
\newtheorem{dfn}[ntn]{Definition}
\theoremstyle{plain}
\newtheorem{lem}[ntn]{Lemma}
\newtheorem{prp}[ntn]{Proposition}
\newtheorem{thm}[ntn]{Theorem}

\theoremstyle{remark}

\DeclareMathAlphabet{\mathds}{U}{dsrom}{m}{n}
\DeclareMathAlphabet{\mathsc}{U}{rsfs}{m}{n}

\newcommand{\pd}[1][]{\partial_{#1}}
\newcommand{\ds}{{\pd[s]}}
\newcommand{\dt}{{\pd[t]}}
\newcommand{\dC}{\mathds{C}}

\newcommand{\dQ}{\mathds{Q}}
\newcommand{\dZ}{\mathds{Z}}
\newcommand{\fs}{{[\![s]\!]}}
\newcommand{\fx}{{[\![\ul x]\!]}}
\newcommand{\fxs}{{[\![\ul x,s]\!]}}
\newcommand{\inv}[1]{{{#1}^{-1}}}
\newcommand{\ms}{{\{\!\{s\}\!\}}}
\newcommand{\ol}[1]{\overline{#1}}
\newcommand{\rd}{\mathrm{d}}
\newcommand{\ri}{\mathrm{i}}
\newcommand{\rH}{\mathrm{H}}
\newcommand{\rR}{\mathrm{R}}
\newcommand{\rM}{\mathrm{M}}
\newcommand{\sC}{\mathsc{C}}

\newcommand{\sG}{\mathsc{G}}
\newcommand{\sH}{\mathsc{H}}
\newcommand{\sL}{\mathsc{L}}
\newcommand{\sO}{\mathsc{O}}
\newcommand{\ul}[1]{{\underline{#1}}}
\newcommand{\wh}[1]{{\widehat{#1}}}
\newcommand{\xymat}{\SelectTips{cm}{}\xymatrix}
\newcommand{\Cd}{{\dC\{\!\{\inv\dt\}\!\}}}
\newcommand{\Cfs}{{\dC\fs}}
\newcommand{\Cfx}{{\dC\fx}}
\newcommand{\Cfxs}{{\dC\fxs}}
\newcommand{\Cs}{{\dC\ms}}
\newcommand{\Ct}{{\dC\{t\}}}
\newcommand{\Cx}{{\dC\{\ul x\}}}

\DeclareMathOperator{\lead}{lead}
\DeclareMathOperator{\ord}{ord}
\DeclareMathOperator{\krn}{ker}
\DeclareMathOperator{\img}{im}
\DeclareMathOperator{\Aut}{Aut}
\DeclareMathOperator{\End}{End}

\newcolumntype{P}[1]{>{\tiny\addtolength{\baselineskip}{2mm}}p{#1}<{\vspace{1mm}}}

\begin{document}

\title{Monodromy of Hypersurface Singularities}

\author{Mathias Schulze}

\address{Mathias Schulze, Department of Mathematics, D-67653 Kaiserslautern}

\email{mschulze@mathematik.uni-kl.de}

\begin{abstract}
We describe algorithmic methods for the Gauss-Manin connection of an isolated hypersurface singularity based on the microlocal structure of the Brieskorn lattice.
They lead to algorithms for computing invariants like the monodromy, the spectrum, the spectral pairs, and M. Saito's matrices $A_0$ and $A_1$.
These algorithms use a normal form algorithm for the Brieskorn lattice, standard basis methods for power series rings, and univariate factorization.
We give a detailed description of the algorithm to compute the monodromy.
\end{abstract}

\maketitle

\tableofcontents

\section{Introduction}

We consider a germ of a holomorphic map $\xymat{f:(\dC^{n+1},0)\ar[r]&(\dC,0)}$ with isolated critical point and Milnor number $\mu$.
J. Milnor \cite{Mil68} first studied this situation by differential geometry.
The regular fibres of a good representative over a punctured disc form a $\sC^\infty$ fibre bundle with fibres of homotopy type of a bouquet of $\mu$ $n$-spheres.
The cohomology of the fibres form a flat vector bundle and there is an associated flat connection on the corresponding sheaf of holomorphic sections, the Gauss-Manin connection.
Moreover, there is a monodromy representation of the fundamental group of the base in the cohomology of the general fibre.
A counterclockwise generator acts via the monodromy which is an automorphism defined over the integers.
\begin{figure}[h]
\caption{The Milnor fibration}
\begin{center}
\epsfig{file=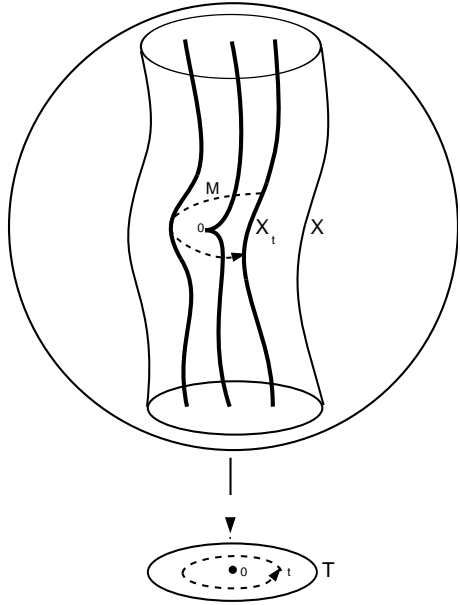,height=8cm}
\end{center}
\end{figure}
Using the De Rham isomorphism, the cohomology of the fibres can be described in terms of holomorphic differential forms.

E. Brieskorn \cite{Bri70} first gave an algorithm to compute the complex monodromy based on this idea.
The original Brieskorn algorithm was first implemented by P.F.M. Nacken \cite{Nac90} in the computer algebra system {\sc Maple V}.
A later implementation \cite{Sch99a,Sch99b} by the author in the computer algebra system {\sc Singular} \cite{GPS01} using standard basis methods turned out to be more powerful.

In section 2 and 3, we briefly introduce the monodromy, the Gauss-Manin connection, the V-filtration, and the Brieskorn lattice, and summarize properties and results which are relevant for the computation of the monodromy.
This will lead us to the notions of saturation and resonnance of lattices in section 4.
In section 5, we describe algorithmic methods for the Gauss-Manin connection based on the microlocal structure of the Brieskorn lattice and the Fourier-Laplace transform \cite{Sch00,SS01}.
We use standard basis methods, univariate factorization, and a normal form algorithm for the microlocal structure of the Brieskorn lattice, the latter of which is not published yet.
These methods lead to algorithms \cite{Sch02,Sch01a,Sch01b} to compute Hodge-theoretic invariants like the spectrum, the spectral pairs, and M. Saito's matrices $A_0$ and $A_1$ \cite{Sai89}.
We describe an algorithm to compute the complex monodromy based on the above ideas.
This is also implemeted in {\sc Singular} \cite{Sch01b} and is much faster then the original Brieskorn algorithm \cite{Sch99b}.
In section 6, we explain how to use the {\sc Singular} implementation by computing an example with a Jordan block of size $3\times3$ of the monodromy.

\section{Monodromy and Gauss-Manin Connection}

Let $\xymat{f:(\dC^{n+1},0)\ar[r]&(\dC,0)}$ be an {\bf isolated hypersurface singularity}.
We choose local coordinates $\ul x=(x_0,\dots,x_n)$ at $0\in\dC^{n+1}$ and $t$ at $0\in\dC$ and set $\ul\pd:=(\pd[0],\dots,\pd[n])$ where $\pd[j]:=\pd[x_j]$.
We denote by
\[
\mu:=\dim_\dC\Cx/\langle\ul\pd(f)\rangle<\infty
\]
the {\bf Milnor number} of $f$.
We choose an $(n+1)$-ball $B$ centered at $0\in\dC^{n+1}$ and a disc $T$ centered at $0\in\dC$ such that
\[
\xymat{B\cap\inv f(0)=:X\ar[r]^-f&T}
\]
is a {\bf Milnor representative} \cite{Mil68}.
Let $\xymat{T':=T\backslash\{0\}\ar@{^(->}[r]^-i&T}$ be the inclusion of the punctured disc.
Then the restriction
\[
\xymat{X\backslash\inv f(0)=:X'\ar[r]^-f&T'}
\]
is a $\sC^\infty$ fibre bundle, the {\bf Milnor fibration}.
The fibres $X_t:=\inv f(T)$ have homotopy type of a bouquet of $\mu$ $n$-sheres and hence the cohomology of the general fibre is given by
\[
\widetilde\rH^k(X_t,\dZ)\cong
\begin{cases}
\dZ^\mu, & k=n,\\
0, & \text{else}.
\end{cases}
\]
The $n$-th {\bf cohomology bundle}
\[
H^n:=\rR^nf_*\dC_{X'}=\bigcup_{t\in T'}\rH^n(X_t,\dC)
\]
is a locally constant sheaf.
The {\bf Gauss-Manin connection} is the associated flat connection
\[
\xymat{\sH^n\ar[r]^-\nabla&\Omega^1_{T'}\otimes_{\sO_{T'}}\sH^n}
\]
on the sheaf of holomorphic sections $\sH^n:=\sO_{T'}\otimes_{\dC_{T'}}H^n$ which is defined by
\[
\nabla(g\otimes v):=\rd g\otimes v,\ g\in\sO_{T'},\ v\in H^n.
\]
We denote by $\xymat{\dt:=\nabla_\dt:\sH^n\ar[r]&\sH^n}$ its covariant derivative with respect to $\dt$.
Lifting paths in $T'$ along flat sections defines an action
\[
\xymat{\pi_1(T',t)\ar[r]&\Aut\bigl(\rH^n(X_t,\dZ)\bigr)}
\]
of the fundamental group $\pi_1(T',t)$ on the $n$-th cohomology of the general fibre.
A counterclockwise generator of $\pi_1(T',t)$ acts by the {\bf monodromy}
\[
\rM\in\Aut\bigl(\rH^n(X_t,\dZ)\bigr).
\]
\begin{thm}[Monodromy Theorem]
The eigenvalues of the monodromy are roots of unity and its Jordan blocks are of size at most $(n+1)\times(n+1)$.
\end{thm}
Let $\xymat{u:T^\infty \ar[r] & T}$, $u(\tau):=\exp(2\pi\ri\tau)$, be the universal covering of $T'$ where $\tau$ is a coordinate on $T^\infty$.
The {\bf canonical Milnor fibre} is defined to be the pullback
\[
X^\infty:=X'\times_{T'}T^\infty
\]
to the universal covering.
The natural maps $\xymat{X_{u(\tau)}\cong X^\infty_\tau\ar@{^(->}[r] & X^\infty}$ are homotopy equivalences.
We consider $A\in\rH^n(X^\infty,\dC)$ as a global flat multivalued section $A(t)$ in $\sH^n$.
Note that
\begin{align*}
\dt A(t)&=0,\\
\rM(A)(\tau)&=A(\tau+1).
\end{align*}

\section{V-Filtration and Brieskorn Lattice}

Let $\rM=\rM_s\rM_u$ be the decompostion of the monodromy into semisimple and unipotent part and set
\[
N:=\log\rM_u.
\]
Note that $N^{n+1}=0$ by the monodromy theorem.
Let
\[
\rH^n(X^\infty,\dC)=\bigoplus_\lambda\rH^n(X^\infty,\dC)_\lambda,\quad\rH^n(X^\infty,\dC)_\lambda:=\ker(\rM_s-\lambda),
\]
be the decomposition into generalized eigenspaces of $\rM$ and
\[
\rM_\lambda:=\rM\arrowvert_{\rH^n(X^\infty,\dC)_\lambda}.
\]
For $A\in\rH^n(X^\infty,\dC)_\lambda$, $\lambda=\exp(-2\pi\ri\alpha)$, $\alpha\in\dQ$, we set
\[
\psi_\alpha(A)(t):=t^\alpha\exp\biggl(-\frac{N}{2\pi\ri}\log t\biggr)A(t).
\]
Then $\psi_\alpha(A)$ is monodromy invariant and hence a global section in $\sH^n$.
\begin{dfn}
We call the $\Ct[\dt]$-module
\[
\sG:=\bigl\langle(i_*\psi_\alpha(A))_0\big\vert\alpha\in\dQ,A\in\rH^n(X^\infty,\dC)_{\exp(-2\pi\ri\alpha)}\bigr\rangle_{\sO_{T,0}}\subset(i_*\sH^n)_0
\]
the {\bf Gauss-Manin connection} of $f$.
\end{dfn}
For all $\alpha\in\dZ$, the map
\[
\xymat{\psi_\alpha:\rH^n(X^\infty,\dC)_\lambda\ar@{^(->}[r]&\sG}
\]
is an inclusion with image $C^\alpha:=\img\psi_\alpha$.
The following lemma shows the correspondence between the monodromy action on $\rH^n(X^\infty,\dC)$ and the $\Ct[\dt]$-module structure on $\sG$ via the maps $\psi_\alpha$.
\begin{lem}\label{1}\
\begin{enumerate}
\item $t\circ\psi_\alpha=\psi_{\alpha+1}$
\item $\dt\circ\psi_\alpha=\psi_{\alpha-1}\circ\bigl(\alpha-\frac{N}{2\pi\ri}\bigr)$
\item $(t\dt-\alpha)\circ\psi_\alpha=\psi_\alpha\circ\bigl(-\frac{N}{2\pi\ri}\bigr)$
\item $\exp(-2\pi\ri t\dt)\circ\psi_\alpha=\psi_\alpha\circ\rM_\lambda$.
\item $C^\alpha=\krn(t\dt-\alpha)^{n+1}$
\item $\xymat{t:C^\alpha\ar[r]&C^{\alpha+1}}$ is bijective.
\item $\xymat{\dt:C^\alpha\ar[r]&C^{\alpha-1}}$ is bijective for $\alpha\ne0$.
\end{enumerate}
\end{lem}
The generalized eigenspaces $C^\alpha$ of the operator $t\dt$ define a filtration on $\sG$.
\begin{dfn}
The {\bf V-filtration} $V$ on $\sG$ is the decreasing filtration by $\Ct$-modules
\[
V^\alpha:=\sum_{\alpha\le\beta}\Ct C^\beta,\quad
V^{>\alpha}:=\sum_{\alpha<\beta}\Ct C^\beta.
\]
\end{dfn}
Note that $V^\alpha/V^{>\alpha}\cong C^\alpha$.
There is not only the $\Ct$-structure on $\sG$.
For $\alpha>-1$, the action of $\inv\dt$ on $V^\alpha$ extends to a structure over a power series ring $\Cd$.
This structure is the key to powerful algorithms.
\begin{dfn}
The ring of {\bf microdifferential operators with constant coefficients} is defined by
\[
\Cd:=\biggl\{\sum_{k\ge0}a_k\pd^{-k}\in\dC[\![\inv\pd]\!]\bigg\vert\sum_{k\ge0}\frac{a_k}{k!}t^k\in\Ct\biggr\}.
\]
\end{dfn}
Note that $\Cd$ is a discrete valuation ring and $t^\alpha\Ct$, $\alpha\in\dQ$, is a free $\Cd$-module of rank $1$.
Together with lemma \ref{1}, this implies the following proposition.
\begin{prp}\ 
\begin{enumerate}
\item For all $\alpha\in\dQ$, $V^\alpha$ is a free $\Ct$-module of rank $\mu$.
\item $\sG$ is a $\mu$-dimensional $\Ct[\inv t]$-vector space.
\item For $\alpha>-1$, $V^\alpha$ is a free $\Cd$-module of rank $\mu$.
\end{enumerate}
\end{prp}
There is a lattice in $\sG$ on which the action of $\dt$ can be computed.
\begin{dfn}
\[
\sH'':=\Omega_{X,0}^{n+1}/\rd f\wedge\rd\Omega_{X,0}^{n-1}.
\]
is called the {\bf Brieskorn lattice}.
\end{dfn}
The Brieskorn lattice is embedded in $\sG$ and is a $\Ct$- and $\Cd$-lattice.
There is an explicit formula for the action of $\dt$ in terms of differential forms.

We summarize these well known properties of the Brieskorn lattice in the following theorem.
\begin{thm}\cite{Seb70,Bri70,Pha77,Mal74}
\begin{enumerate}
\item $\sH''$ is a free $\Ct$-module of rank $\mu$.
\item $\xymat{s:\sH''\ar@{^(->}[r]&\sG}$, $s\bigl([\omega]\bigr)(t):=\int\frac{\omega}{\rd f}\big\vert_{X_t}$
\item $\dt s\bigl([\rd f\wedge\eta]\bigr)=s\bigl([\rd\eta]\bigr)$
\item $\sH''$ is a free $\Cd$-module of rank $\mu$.
\item $V^{-1}\supset\sH''\supset V^{n-1}$
\end{enumerate}
\end{thm}
The action of $\dt$ may have a pole of order up to $n+1$ on $\sH''$.
Since the monodromy is related to the action of $t\dt$, we consider $t\dt$-invariant lattices in the next section.

\section{Saturation and Resonance}

Since $\Ct$ is a discrete valuation ring, for any two $\Ct$-lattices $\sL,\sL'\subset\sG$, there is a $k\in\dZ$ such that $t^k\sL\subset\sL'$.
Hence, for any $\Ct$-lattice $\sL$,
\[
V^{\alpha_2}\subset\sL\subset V^{>\alpha_1}
\]
for some $\alpha_1,\alpha_2\in\dQ$.
\begin{dfn}
Let $\sL\subset\sG$ be a $\Ct$-lattice.
If $t\dt\sL\subset\sL$ then $\sL$ is called {\bf saturated} and the induced endomorphism $\ol{t\dt}\in\End_\dC(\sL/t\sL)$ is called the {\bf residue} of $\sL$.
If $\ol{t\dt}$ has non-zero integer differences of eigenvalues then $\sL$ is called {\bf resonant}.
\end{dfn}

Let $\sL\subset\sG$ be a $\Ct$-lattice with $V^{\alpha_2}\subset\sL\subset V^{>\alpha_1}$ for some $\alpha_1,\alpha_2\in\dQ$.
By the Leibnitz rule and since $V^{\alpha_1}$ is saturated,
\[
\sL_0:=\sL,\quad\sL_{k+1}:=\sL_k+t\dt\sL_k,
\]
defines an increasing sequence of $\Ct$-lattices
\[
V^{\alpha_2}\subset\sL_0\subset\sL_1\subset\cdots\subset V^{\alpha_1}.
\]
Since $V^{\alpha_1}$ is noetherian, this sequence is stationary and
\[
\sL_\infty:=\bigcup_{k\ge0}\sL_k.
\]
is a saturated $\Ct$-lattice.
\begin{dfn}
$\sL_\infty$ is called the {\bf saturation} of $\sL$.
\end{dfn}
The following proposition is not difficult to prove using lemma \ref{1}.
\begin{prp}\cite{GL73}
$\sL_{\mu-1}=\sL_\infty$
\end{prp}
If $\sL$ is saturated then
\[
\sL=\Bigl(\bigoplus_{\alpha_1<\alpha<\alpha_2}\sL\cap C^\alpha\Bigr)\oplus V^{\alpha_2}.
\]
Together with lemma \ref{1}, this implies the following proposition.
\begin{prp}
Let $\sL\subset\sG$ be a saturated $\Ct$-lattice with residue $\ol{t\dt}\in\End_\dC(\sL/t\sL)$.
\begin{enumerate}
\item The eigenvalues of $\exp(-2\pi\ri\ol{t\dt})$ are the eigenvalues of the monodromy.
\item If $\sL$ is non-resonant then $\exp(-2\pi\ri\ol{t\dt})$ is conjugate to the monodromy.
\end{enumerate}
\end{prp}

\section{Microlocal Structure and Algorithms}

We abbreviate $\Omega^\bullet:=\Omega_{X,0}^\bullet$.
The microlocal structure extends to a $\Cs[\ds]$-module structure by the {\bf Fourier-Laplace transform}
\[
s:=\inv\dt,\quad\ds:=\dt^2t=s^{-2}t.
\]
Since\[
[\ds,s]=[\dt^2t,\inv\dt]=\dt^2t\inv\dt-\dt t=1,
\]
$V^\alpha$ for $\alpha>-1$ and $\sH''$ are $\Cs[s^2\ds]$-modules.
Note that $t=s^2\pd[s]$ is a $\Cs$-derivation.
Since
\[
\dt t=\inv st=s\ds,
\]
the saturation $\sL_\infty$ of a $\Ct$- and $\Cs$-lattice $\sL$ is a saturated $\Ct$- and $\Cs$-lattice.
Note that this holds for $\sH''$.

By the finite determinacy theorem, we may assume that $f\in\dC[x]$ is a polynomial.
A $\dC$-basis of
\[
\sH''/s\sH''=\Omega^{n+1}/\rd f\wedge\Omega^n\cong\Cx/\langle\ul\pd(f)\rangle
\]
represents a $\Cs$-basis of $\sH''$.
If $\ul g$ is a standard basis of $\langle\ul\pd(f)\rangle$ with respect to a local monomial ordering then the monomials which are not contained in the leading ideal $\langle\lead\ul g\rangle=\langle\lead(\ul\pd(f))\rangle$ form a monomial $\dC$-basis of $\Cx/\langle\ul\pd f\rangle$.
Hence, one can compute a monomial $\Cs$-basis $\ol m=\begin{pmatrix}m^1\\\vdots\\m^\mu\end{pmatrix}$ of $\sH''$.

We define the {\bf $\ol m$-matrix $A=A(s)=\sum_{k\ge 0}A_ks^k$ of $t$} by
\[
t\ol m=:A\ol m.
\]
Then the $\ol m$-basis representation of $t$ on $\sH''$ is given by
\[
t\ul g\ol m=\bigl(\ul gA+s^2\ds(\ul g)\bigr)\ol m.
\]
If $U$ is a $\Cs[\inv s]$-basis transformation and $A'$ the $U\ol m$-matrix of $t$ then
\[
A'=\bigl(UA+s^2\ds(U)\bigr)\inv U
\]
is the basis transformation formula for $U$.

Let $\wh\Omega$ resp. $\wh\sH''$ be the $\langle\ul x\rangle$-adic resp. $\langle s\rangle$-adic completion of $\Omega$ resp. $\sH''$.
The isomorphism $\xymat{\rd\ul x:\Cfx\ar[r]_-\sim&\wh\Omega^{n+1}}$ induces an isomorphism
\[
\Cfxs\bigg/\sum_{j=0}^n\bigl(\pd[j](f)-s\pd[j]\bigr)\Cfxs\cong_\Cfs\wh\sH''
\]
and $\wh\sH''$ is a {\bf differential deformation} of the Jacobian algebra $\Cfx/\langle\ul\pd(f)\rangle$ in the sense of \cite{Sch01c}.
Using the {\bf normal form algorithm} in \cite{Sch01c}, one can compute any $K$-jet $A_{\le K}=A_{\le K}(s)=\sum_{k=0}^KA_ks^k$ of $A$.

The $\ol m$-basis representation $H''_\infty$ of the saturation $\sH''_\infty$ can be computed recursively by
\begin{align*}
H''_0&:=\Cs^\mu,\\
H''_{k+1}&:=H''_k+\bigl(\inv sH''_kA_{\le k}+s\ds H''_k\bigr).
\end{align*}
Note that only finite jets of $A$ are involved.
We use a local monomial degree ordering.
By computing a normal form of $H''_{k+1}$ with respect to $H''_k$, one can check when the sequence
\[
H''_0\subset H''_1\subset H''_2\subset\cdots
\]
becomes stationary and compute generators of the saturation $H''_\infty$ of $H''_0$.

By Nakayama's lemma, a minimal standard basis $M=\begin{pmatrix}\ul m^1\\\vdots\\\ul m^\mu\end{pmatrix}$ of $H''_\infty$ is a $\Cs$-basis.
The {\bf $M\ol m$-matrix $A'=A'(s)=\sum_{k\ge0}A'_ks^k$ of $t$} is defined by
\[
MA+s^2\ds M=:A'M.
\]
We set
\[
\delta(M):=\max\bigl\{\ord\bigl(m^i_j\bigr)-\ord\bigl(m^k_l\bigr)\big\vert m^i_j\ne0\ne m^k_l\bigr\}\le\mu-1
\]
such that
\[
\big(MA_{\le K+\delta(M)}+s^2\ds M\big)_{\le K}=:A'_{\le K}M.
\]
for any $K\ge0$.
Note that only finite jets of $A$ are involved.
Hence, one can compute any $K$-jet $A'_{\le K}=A'_{\le K}(s)=\sum_{k=0}^KA'_ks^k$ of $A'$.
Note that the $M\ol m$-basis representation of $t$ on $\sH''_\infty$ is given by
\[
t\ul gM\ol m=\bigl(\ul gA'+s^2\ds(\ul g)\bigr)M\ol m.
\]
Hence, the $M\ol m$-basis representation of $\dt t=\inv st$ on $\sH''_\infty/s^K\sH''_\infty$ is given by
\[
\dt t\ul gM\ol m=\bigl(\inv s\ul gA'_{\le K}+s\ds(\ul g)\bigr)M\ol m
\]
and $\inv sA'_{\le 1}=A'_1$ is the $M\ol m$-basis representation of the residue of $\sH''_\infty$.
Note that the eigenvalues of $A'_1$ are rational by the monodromy theorem and can be computed using univariate factorisation.
If $A'_1$ is non-resonant then $\exp(-2\pi\ri A'_1)$ is a monodromy matrix.

Otherwise, we proceed as follows.
Let $\delta(A')>0$ be the maximal integer difference of $A'_1$.
First we compute $A'_{1+\delta(A')}$ from $A_{1+\delta(M)+\delta(A')}$ as before.
After a $\dC$-linear coordinate transformation, we may assume that
\[
A'=\begin{pmatrix}{A'}^{1,1}&{A'}^{1,2}\\{A'}^{2,1}&{A'}^{2,2}\end{pmatrix}
\]
with ${A'}^{1,2}_1=0$, ${A'}^{2,1}_1=0$, $A'_0=0$, the eigenvalues of ${A'}^{1,1}$ are minimal in their class modulo $\dZ$, and the eigenvalues of ${A'}^{2,2}$ are non-minimal in their class modulo $\dZ$.
Then the $\Cs[\inv s]$-coordinate transformation
\[
U:=\begin{pmatrix}s&0\\0&1\end{pmatrix}
\]
gives
\[
A''=\begin{pmatrix}{A''}^{1,1}&{A''}^{1,2}\\{A''}^{2,1}&{A''}^{2,2}\end{pmatrix}
:=\bigl(UA'+s^2\ds(U)\bigr)\inv U
=\begin{pmatrix}{A'}^{1,1}+s&s{A'}^{1,2}\\\inv s{A'}^{2,1}&{A'}^{2,2}\end{pmatrix}.
\]
Note that $A''_0=0$, $\delta(A'')\le\delta(A')-1$, and that $A''_{K}$ depends only on $A'_{\le K+1}$.
With $\delta(A')$ of these transformations we decrease $\delta(A')$ to zero such that $\exp(-2\pi\ri A'_1)$ is a monodromy matrix as before.

The above methods lead to algorithms \cite{Sch02,Sch01a,Sch01b} to compute Hodge-theoretic invariants like the spectrum, the spectral pairs, and M. Saito's matrices $A_0$ and $A_1$ \cite{Sai89}.

\section{Examples}

The algorithms described above is implemented in the computer algebra system {\sc Singular} \cite{GPS01} in the library {\tt gaussman.lib} \cite{Sch01a}.
We use this implementation to compute an example.

First, we have to load the library:
\begin{verbatim}
> LIB "gaussman.lib";
\end{verbatim}
Then we define the ring $R:=\dQ[x,y,z]_{(x,y,z)}$ and the polynomial $f=x^2y^2z^2+x^7+y^7+z^7\in R$:
\begin{verbatim}
> ring R=0,(x,y,z),ds;
> poly f=x2y2z2+x7+y7+z7;
\end{verbatim}
Finally, we compute the Jordan data of the monodromy of the singularity defined by $f$ at the origin.
\begin{verbatim}
> spprint(monodromy(f));
((1/2,1),18),((1/2,3),1),((9/14,1),15),((9/14,2),3),
((11/14,1),15),((11/14,2),3),((6/7,1),3),((13/14,1),15),
((13/14,2),3),((1,2),1),((15/14,1),15),((15/14,2),3),
((8/7,1),3),((17/14,1),15),((17/14,2),3),((9/7,1),3),
((19/14,1),15),((19/14,2),3),((10/7,1),3),((11/7,1),3),
((12/7,1),3)
\end{verbatim}
The computation takes about $2$ minutes in a {\sc Pentium\,III\,800}.
A Jordan block of the monodromy of size $s$ with eigenvalue $\exp(-2\pi\ri\alpha)$ occuring with multiplicity $m$ is denoted by $((\alpha,s),m)$.
Note that there is a Jordan block of size $3$ with eigenvalue $-1$ which is the maximum possible size.

\nocite{AGV88}
\nocite{Nac90}
\nocite{Pha79}
\nocite{Sai89}
\nocite{Sch96}
\nocite{SS85}
\bibliographystyle{amsalpha}
\bibliography{mohs}

\end{document}